\documentclass[12pt]{amsart}
\usepackage{amssymb}
\theoremstyle{plain}
 \newtheorem{thm}{Theorem}[section]

\theoremstyle{definition}

\theoremstyle{remark}
 \newtheorem{rem}{Remark}[section]
\begin{document}
\title[A note on 
rigidity for crossed 
product 
von Neumann algebras]
{A note on 
rigidity for crossed 
product 
von Neumann algebras}
\author[Tomohiro Hayashi]{{ Tomohiro Hayashi} }
\address{Nagoya Institute of Technology, 
Gokiso-cho, Showa-ku, Nagoya, Aichi, 466-8555, Japan}
\email{hayashi.tomohiro@nitech.ac.jp}
\thanks{The author was supported by 
Grant-in-Aid for Young Scientists 
(B) of Japan Society for 
the Promotion of Science.}
\baselineskip=17pt
\maketitle

\begin{abstract} 
In this note, we will point 
out, 
as a corollary of Popa's rigidity 
theory, 
that 
the crossed product von Neumann algebras 
for Bernoulli shifts cannot have 
relative property T. 
This is an 
operator algebra analogue of the 
theorem shown by Neuhauser and 
Cherix-Martin-Valette for discrete 
groups. Our proof is 
different from that for groups. 
\end{abstract} 

\section{Introduction}
Nowadays, the notion of relative property T plays 
one of the 
central roles for operator algebras 
(See \cite{PP} for information). 
The most cerebrated results in this area is 
the rigidity theory for malleable 
actions, established by 
Popa \cite{P3}\cite{P4}. For example, Popa solved Connes' problem 
for a certain class of discrete groups by using his 
rigidity theory. 
In this note, we would like to 
point out that as a corollary of Popa's rigidity 
theory, one can show the operator algebra analogue of the 
theorem shown by Neuhauser and 
Cherix-Martin-Valette for discrete 
groups. 

Recently Neuhauser and 
Cherix-Martin-Valette 
independently showed the following 
theorem: (See \cite{BHV}\cite{CMV}
\cite{N} for the proof.)

\begin{thm}[Neuhauser, 
Cherix-Martin-Valette]
Let $G$, $H$ be countable discrete groups with 
$H\not=\{e\}$. If $G$ is an 
infinite group, then the 
inclusion 
$\oplus_{G}H\subset
(\oplus_{G}H)\rtimes G$ 
does not have 
relative property T.
\end{thm} 
\noindent Then it is natural to consider the 
operator algebra analogue 
of this result. The main result of this paper 
is: 

\begin{thm} 
Let $A(\not={\Bbb C})$ 
be a finite von Neumann algebra with a distinguished 
faithful normal tracial state $\tau$. Let $G$ be a countable 
discrete, infinite group. 
Consider the finite von Neumann algebra 
$Q=\otimes_{G}A$ with a trace $\otimes_{G}\tau$. Let 
$\sigma$ be the Bernoulli shift action of $G$ on $Q$. 
Then the inclusion 
$Q\subset P=Q\rtimes_{\sigma}G$ 
does not have relative property T 
in the sense of \cite{PP}. 

Moreover, if $A$ is type I, 
for any diffuse von Neumann subalgebra 
$Q_{0}\subset Q$, the inclusion 
$Q_{0}\subset P=Q\rtimes_{\sigma}G$ 
does not have relative property T. 
\end{thm} 
Here we remark that if $A$ itself is a property T 
$II_{1}$ factor for example, 
then obviously the inclusion 
$A\subset P$ 
has relative property T. This means that in the above theorem, 
we cannot omit the condition that $A$ 
is type I.

The proof for the group case is elegant 
and understandable 
\cite{BHV}\cite{CMV}. However 
it is difficult (at least for the author) to 
arrange it for the operator algebra 
setting. Thus we shall adopt a slightly 
different method. Our main tool is Popa's 
rigidity theory, as noted above. 

\section{Proof}
Before starting the proof, we shall 
fix some notations. 
For a finite von Neumann algebra 
$M$ with a faithful normal 
tracial state $\tau$, we denote 
its 2-norm by 
$||x||_{2}=\tau(x^{*}x)^{1/2}$ for 
$x\in M$. The norm-unit ball of $M$ 
is denoted by $(M)_{1}$. 
That is, $(M)_{1}=\{x\in M;\ 
||x||\leq 1\}$. For a von Neumann 
subalgebra $N\subset M$, we denote 
by $E_{N}$ 
the trace-preserving conditional expectation 
onto $N$. The center of $M$ is denoted by 
$Z(M)$. For a discrete group $G$, 
we denote its group von Neumann algebra by 
$L(G)$. If $G$ acts on $M$, then $L(G)$ can be 
considered as a subalgebra of the 
crossed product von Neumann algebra 
$M\rtimes G$ in the natural way. 
We denote the canonical implementing 
unitary in $L(G)$ by 
$\lambda_{g}$ for 
any $g\in G$. 

Let us now start the proof of 
theorem 1.2. 

First we would like to give the proof 
that $Q\subset P$ does not have 
relative property T. This can be very much easily 
shown and the proof might be a folklore 
for specialists. I learned the following 
proof from the referee. I would like to 
thank him. 

For any positive number $r>0$, 
let $\phi_{r}(\cdot)=r\tau(\cdot)+
(1-r){\rm{id}}$ and let 
$\Phi_{r}(\cdot)=\otimes_{G}\phi_{r}(\cdot)$. 
Then it is easy to see that this unital, normal, 
trace-preserving 
completely positive map on $Q$ can be extended to 
$P$ by 
$\Phi_{r}(x\lambda_{g})=\phi_{r}(x)\lambda_{g}$ 
for any $x\in Q$ and $g\in G$. We use the same notation 
$\Phi_{r}$ for this extended map. 
If the inclusion $Q\subset P$ 
is rigid, then $\Phi_{r}$ must converge 
to identity uniformly 
on $(Q)_{1}$ with respect to 
the trace 2-norm $||\cdot||_{2}$ 
as $r\rightarrow 0$. 
But this is impossible. Indeed, 
Let $u\in A$ be some non-scalar unitary element 
so that $|\tau(u)|<1$. 
(Since $A$ is not scalar, we can always find such a unitary.) 
Then
$\|\phi_r(u)\|_{2}=\|r\tau(u)+(1-r)u\|_{2}\leq 
r|\tau(u)|+1-r$. 
For any finite subset $F\subset G$, 
define $u_{F}=\otimes_{F}u\in Q$. 
Then we have 
$\|\Phi_r(u_F)\|_2\leq
(r|\tau(u)|+1-r)^{|F|}\leq 1$. 
Since 
$\lim_{r\rightarrow 0}\|\Phi_r(u_F)\|_2=
\|u_{F}\|_{2}=1$ uniformly for the choice of 
$F$, we see that 
$\lim_{r\rightarrow 0}(r|\tau(u)|+1-r)^{|F|}
=1$ uniformly for the choice of 
$F$. This is obviously impossible.

Next we shall consider the case that 
$A$ is type I. Let $Q_{0}$ 
be a diffuse von Neumann subalgebra 
of $Q=\otimes_{G}A$. We show that 
$Q_{0}\subset P$ does not have 
relative property T. 
The following theorem plays a crucial role 
in our argument: 

\begin{thm}[Popa \cite{P3}, 
special case (uniqueness of group algebra)] 
Let $C=A\otimes B$ where 
$A$ is a type I von Neumann algebra 
with a faithful normal tracial state 
$\tau_{1}$ and 
$B$ is diffuse abelian with 
faithful normal tracial state 
$\tau_{2}$. 
Let $G$ be a countable infinite 
group. 
($G$ may be non-property T!) 
Consider the Bernoulli shift action 
$\sigma$ of $G$ on $N=\otimes_{G}C$, 
where the von Neumann algebra 
$N=\otimes_{G}C$ is a 
completion 
with respect to $\tau=\otimes_{G}(\tau_{1}\otimes\tau_{2})$. 
Let $P\subset N\rtimes_{\sigma}G$ be an 
irreducible $II_{1}$-subfactor. 
Under these assumptions, if there exists a diffuse 
von Neumann subalgebras $Q\subset P$ 
such that this inclusion is rigid, 
(for example, if $P$ has property T, 
we can take $Q=P$,) 
then there exists a non-zero partial 
isometry $v\in N\rtimes_{\sigma}G$ 
such that 
$vQv^{*}\subset L(G)$. 
\end{thm}

We shall give the proof of theorem 1.2. 
in the case that 
$A$ is type I. 

Let $A$, $G$, $Q$, $Q_{0}$, $P$, $\sigma$ be 
as in theorem 1.2 and we assume 
$A$ is type I. 
Let $B=L^{\infty}[0,1]$ and $C=A\otimes B$. 
We identify 
$P$ with 
$(\otimes_{G}A\otimes 1)\rtimes G 
\subset M=\otimes_{G}C\rtimes G$. Then 
by assumption, 
$P$ is an irreducible 
$II_{1}$-subfactor of $M$. 
(Irreducibility follows from 
$L(G)'\cap M=Z(L(G))\subset P$.) 
Suppose that $Q_{0}\subset P$ 
has relative property T. 
Since $Q_{0}$ is diffuse, 
by Popa's theorem we can find a non-zero 
partial isometry $v\in M$ such that 
$vQ_{0}v^{*}\subset L(G)$. But this is 
obviously impossible. 
Indeed it is easily seen, for example, 
by using 
Popa's orthogonal pair technique~\cite{P2} 
as follows: 
Since $Q_{0}\subset\otimes_{G}A$ is diffuse, 
for any $\epsilon>0$, we can take a family 
of orthogonal projections 
$\{e_{i}\}_{i=1}^{n}\subset Q_{0}$ such that 
$\sum_{i}e_{i}=1$ 
and $\tau(e_{i})<\epsilon$. 
Then, for any $x\in \otimes_{G}C$ and 
$g\in G$ we see that 
\begin{align*} 
|\tau(v^{*}x\lambda_{g})|^{2}&=
|\tau(\sum_{i}e_{i}v^{*}x\lambda_{g}e_{i})|^{2}\\
&\leq 
||\sum_{i}e_{i}v^{*}x\lambda_{g}e_{i}||_{2}^{2}\\
&=\sum_{i}||e_{i}v^{*}x\lambda_{g}e_{i}||_{2}^{2}\\
&=\sum_{i}\tau(ve_{i}v^{*}x\sigma_{g}(e_{i})x^{*})\\
&=\sum_{i}\tau(ve_{i}v^{*})\tau(x\sigma_{g}(e_{i})x^{*})\\
&\leq \epsilon ||x||_{2}. 
\end{align*} 
Here we use the fact that 
$ve_{i}v^{*}\in L(G)$ and 
$x\sigma_{g}(e_{i})x^{*}\in \otimes_{G}C$. 
Since $\epsilon$ is arbitrary, $v$ must be $0$. 
This is a contradiction and hence 
the inclusion $Q_{0}\subset P$ 
does not have relative property T. 

\begin{rem} 
The proof in the group case (Theorem 1.1) is 
done by constructing the non-vanishing cocycle 
of $G$ to some Hilbert space explicitly. 
It is difficult (for the author) to modify 
this argument for the von Neumann algebra case. 
\end{rem}

\end{document}